# SPECTRAL DYNAMICS OF THE VELOCITY GRADIENT FIELD IN RESTRICTED FLOWS

HAILIANG LIU AND EITAN TADMOR

ABSTRACT. We study the velocity gradients of the fundamental Eulerian equation, $\partial_t u + u \cdot \nabla u = F$, which shows up in different contexts dictated by the different modeling of $F$'s. To this end we utilize a basic description for the spectral dynamics of $\nabla u$, expressed in terms of the (possibly complex) eigenvalues, $\lambda = \lambda(\nabla u)$, which are governed by the Ricatti-like equation $\lambda_t + u \cdot \nabla \lambda + \lambda^2 = \langle l, \nabla Fr \rangle$.

We focus our investigation on four prototype models associated with different forcing $F$, ranging from simple linear damping and a viscous dusty medium models to the main thrust of the paper — the restricted models of Euler/Navier-Stokes equations and Euler-Poisson equations.

In particular, we address the question of the time regularity for these models, that is, whether they admit a finite time breakdown, a global smooth solution, or an intermediate scenario of critical threshold phenomena where global regularity depends on initial configurations.

Using the spectral dynamics as our essential tool in these investigations, we obtain a simple form of a critical threshold for the linear damping model and we identify the 2D vanishing viscosity limit for the viscous irrotational dusty medium model. Moreover, for the $n$-dimensional restricted Euler equations we obtain $[n/2] + 1$ global invariants, interesting for their own sake, which enable us to precisely characterize the local topology at breakdown time, extending previous studies in the $n = 3$-dimensional case. Finally, as a forth model we introduce the $n$-dimensional restricted Euler-Poisson (REP)system, identifying a set of $[n/2]$ global invariants, which in turn yield (i) sufficient conditions for finite time breakdown, and (ii) characterization of a large class of 2-dimensional initial configurations leading to global smooth solutions. Consequently, the 2D restricted Euler-Poisson equations are shown to admit a critical threshold.

## Contents

1. INTRODUCTION

It is well known that the velocity gradients in a turbulent flow are larger than their mean gradients by at least a factor of order $\sqrt{R_\delta}$, with $R_\delta$ being the Reynolds number based on internal length and velocity scales. Fluctuation gradients are limited by the mean flow and contribute a dominant portion of the kinetic energy dissipation, but otherwise they contribute nothing to the mean transport of momentum because of the linearity of the viscous stress term in the Navier-Stokes equations. Consequently, much research has been directed at gaining a better understanding of the velocity gradient field, $\nabla u$, which is completely dictated by the vorticity in incompressible flows, [2, 9, 19, 30, 6].

Motivated by such questions, the goal of this work is to present new observations on the velocity gradients for a general class of so called restricted flows, where the velocity field, $u$, is governed by the Newtonian law,

$$(1.1) \qquad \partial_t u + u \cdot \nabla u = F, \quad x \in \mathsf{IR}^n,$$

with $F$ being a general forcing acting on the flow. Different regimes of the flow are modeled by different $F's$. A key issue in this line of research is the control of the velocity gradient $\nabla u$, and a classical approach in this context, is to consider linear combinations of the entries of $\nabla u$, controlling physically relevant quantities like vorticity, divergence, etc, see [2, 10, 30].

The novelty of the analysis taken in the present article is the use of the eigenvalues of the velocity gradient field. The eigenvalues, $\lambda = \lambda(\nabla u)$, exhibit of course a strong nonlinear dependence on the entries of $\nabla u$, and are shown to play a crucial role in governing the behavior of the flow. Indeed, the dynamics of these eigenvalues $\lambda(\nabla u)$, is shown, in §3, to be governed by Ricatti-like equation

$$\partial_t \lambda + u \cdot \nabla \lambda + \lambda^2 = <l, \nabla Fr>,$$

with $l(r)$ being the left (right) eigenvectors of $\nabla u$. Equipped with this description for the spectral dynamics of $\nabla u$, we turn to study several physical models with different forcing, outlined in §2 and analyzed in §4–§7.

We focus on four prototype models in this paper. The first two are a simple linear damping model studied in §4, and a viscous dusty medium model in §5. Next, the main thrust of the paper is devoted to the restricted models of Euler equations in §6, and in §7 we introduce the so called restricted Euler-Poisson equations as our fourth model problem. We focus our attention on the question of time regularity for these models, that is, whether they admit a finite time breakdown, a global smooth solution or an intermediate scenario of critical threshold phenomena where global regularity depends on the initial configurations as in e.g. [15].

The question of time regularity is of fundamental importance from both mathematical and physical points of view, and a considerable effort is still being devoted to this issue for both compressible and incompressible Euler equations. Consult [31, 20, 18, 2, 9, 33, 26] for a partial list of recent references. In particular, the possible phenomena of finite time breakdown for 3D incompressible flows signifies the onset of turbulence in higher Reynolds number flows. Several simplified models for 3-D Euler equations were proposed to understand this phenomena, see [34] for a restricted dynamics model, [10] for a vorticity dynamics model, [11] for a so called distorted model as well as a stochastic model in [12].

The paper is organized as follows.

After introducing the basic spectral dynamics Lemma 3.1 in §3, we begin our discussion of the time regularity with the simple linear damping forcing model in §4.

In §5 we deal with viscosity forces, where we study the irrotational viscous dusty medium model, and identify its 2D vanishing viscosity limit. Here, the spectral dynamics offers us a novel approach at the level of the velocity field, $u = \nabla \phi$ — an alternative to the classical notion of viscosity solutions for Hamilton-Jacobi equations at the level of the potential $\phi$. Spectral dynamics serves as an essential tool in our approach, most notably the use of a key a priori $L^1$-contraction estimate expressed in terms of the unintuitive nonlinear quantity $\lambda_2 - \lambda_1$, with $\lambda_i = \lambda_i(\nabla u), i = 1, 2$ being the two real eigenvalues of $\nabla u$.

In §6 we use the spectral dynamics to revisit the restricted Euler models introduced by Vieillefosse in [34]. The so called restricted Euler equations (RE for short) refer to a localized model of the Euler/Navier-Stokes equations, where the usual global pressure forces are replaced by their local, isotropic trace. We study the time regularity of the general $n$-dimensional RE equations, extending the previous studies in [34, 5] for the special $n = 3$ case. Here, we enjoy the advantage of using the spectral dynamics of such $n$-dimensional flows, which enables us to identify a large set of at least $[n/2] + 1$ independent integrals of the motion. Using these $[n/2]+1$ global invariants, interesting for their own sake, we precisely characterize the finite time breakdown for the $n$-dimensional RE equations.

We note in passing that the RE model has been an appealing candidate for describing the dynamics of the local velocity gradient, [1, 3]. Despite this restricted approximation to the pressure, the RE equations can still describe the local topology of Euler equations and capture certain statistical features of the physical flow. In this spirit we introduce, in §7, a restricted model for the Euler-Poisson system, so called REP equations. For general $n$-dimensional REP equations we obtain a set of at least $[n/2]$ global invariants, which in turn yields

(1) Sufficient conditions for the finite time breakdown in $n$-dimensional REP equations. Moreover, we characterize the precise local topology of the flow at breakdown time;

(2) Sufficient conditions for a large class of 2D initial configurations leading to the existence of global smooth solution for 2D REP equations.

We point out that though RE model was sometimes argued for its unphysical finite time singularity, our REP model does support the global smooth solutions. In particular, it follows that the 2D REP equations admit a critical threshold, distinguished between initial configurations leading to either the finite time breakdown or the global smooth solutions. We refer to [29] for a detailed study of such phenomena for this 2D REP model.

In §8 we discuss possible extensions of the results obtained in this work and we comment on some remaining open issues. Finally, in the Appendix we revisit the spectral dynamics of the general $n \times n$ RE models from yet another perspective of a trace dynamics, extending the study of traces, $tr(\nabla u)^k$, $k = 1, \cdots, n$, initiated in [34] for the special $n = 3$ case.

## 2. Basic equations — four prototype models

In what follows we shall require the equations governing the dynamics of a fluid in both the Eulerian and Lagrangian forms. We shall study the flow of a fluid which initially at $t = 0$ occupies the whole space $\mathrm{I\!R}^n$ for arbitrary dimension $n$, although only the cases $n = 2$ and $n = 3$ have a clear physical meaning.

Let a Cartesian coordinate system be fixed in $\mathrm{I\!R}^n$. We denote by $\alpha$ the initial position of a fluid particle. The motion of the fluid is assumed to be given, if for any $\alpha \in \mathrm{I\!R}^n$ the position $x(\alpha, t) \in \mathrm{I\!R}^n$ of the fluid particle is known for all $\alpha \in \mathrm{I\!R}^n$ and for all $t \in \mathrm{I\!R}^+$, with $x(\alpha, 0) = \alpha$. Further, $\alpha \to \frac{d}{dt} x(\alpha, t) = u(x, t)$ is a Lagrangian velocity field at the time $t$. The Lagrangian equations of the dynamics of a fluid amount to

$$\frac{d^2}{dt^2} x = F,$$

where $F$ denotes the forcing acting on the fluid. The corresponding Eulerian equations in the standard form read

(2.1) $$\partial_t u + u \cdot \nabla u = F, \quad x \in \mathrm{I\!R}^n,$$

where $\frac{d}{dt} = \frac{\partial}{\partial t} + u \cdot \nabla$ is the Lagrangian derivative. Equation (2.1) shows up in a variety of contexts dictated by the different modeling of $F$'s.

Differentiation of the above equation with respect to $x$ gives the relation for the local velocity gradient tensor, $M := \nabla u$,

(2.2) $$\partial_t M + (u \cdot \nabla) M + M^2 = \nabla F.$$

The central issue of interest here is to control the local velocity gradient tensor in (2.2) and to clarify whether the associated distortion matrix, $\Gamma := \partial x / \partial \alpha$, remains nonsingular as time evolves. In particular, one is interested to know whether there is a finite time breakdown, a global smooth solution or an intermediate scenario of critical threshold phenomena with a conditional breakdown, consult e.g.,[15]. In the remainder of this section we discuss four prototype examples associated with different forcing $F$. In the following sections, we will follow the spectral dynamics of the velocity gradient tensor associated with these four examples to demonstrate the different phenomena of global regularity, finite time breakdown as well as the existence of critical threshold.

### 2.1. Linear damping.
Consider a model of the form

(2.3) $$\partial_t u + u \cdot \nabla u = Cu, \quad x \in \mathrm{I\!R}^n.$$

Here we deal with the simple forcing, $F = Cu$, where $C$ is a constant matrix. The corresponding local velocity tensor $M = \nabla u$ solves

(2.4) $$\partial_t M + (u \cdot \nabla) M + M^2 = CM.$$

In §4 we use the spectral dynamics of $M$ to show that there exists a critical threshold depending on the choice of the matrix $C$.

## 2.2. Irrotational viscous flow.

Next we consider the viscous forces, $F := \nu \Delta u$, which leads to the so called viscous dusty medium model

$$\partial_t u + u \cdot \nabla u = \nu \Delta u, \quad x \in \mathbb{R}^n, \tag{2.5}$$

where $\nu > 0$ is a viscosity amplitude. Other suggested names are Burgers system [16], Hopf system, Riemann equation (for $n = 1$). Zeldovitch [36] proposed to consider this system as a model describing the evolution of the rarefied gas of non-interacting particles.

The Hopf-Cole transformation, $u = -2\nu \nabla[log(\psi)]$, links the Burgers system to the heat equation

$$\partial_t \psi = \nu \Delta \psi$$

provided the initial data, $u_0 = u(x,0)$, admits the form $u_0 = -2\nu\nabla[log(\psi_0)]$ for some positive $\psi_0 = \psi(x,0)$ (this is available for $n = 1$). The corresponding local velocity gradient field satisfies

$$\partial_t M + (u \cdot \nabla)M + M^2 = \nu \Delta M, \quad x \in \mathbb{R}^n. \tag{2.6}$$

We focus our attention on solutions to the 2D irrotational viscous flows, $u = u^\nu$, and we use the spectral dynamics of $M$ in order to study the inviscid limit, $u = \lim_{\nu \to 0} u^\nu$. In particular, the limiting 2D irrotational flow is shown to be a weak solution of

$$\partial_t u + u \cdot \nabla u = 0, \quad u =: \nabla \phi,$$

which is interpreted through the Eikonal equation $\partial_t \phi + |\nabla \phi|^2/2 = 0$.

## 2.3. Restricted Euler/Navier Stokes equations.

For the forcing involving viscosity and pressure, we consider the Navier-Stokes equations of incompressible fluid flow in $n$ space dimensions, which can be expressed as the system of $n+1$ equations,

$$\partial_t u + u \cdot \nabla u = \nu \Delta u - \nabla p, \quad x \in \mathbb{R}^n, \quad t > 0, \tag{2.7}$$

$$\nabla \cdot u = 0, \tag{2.8}$$

$$u(x,0) = u_0(x). \tag{2.9}$$

Here $u$ is the fluid velocity, $p$ is the scalar pressure, and $\nu > 0$ is the reciprocal of the Reynolds number. When the coefficient $\nu$ vanishes in (2.7), we have the incompressible Euler equations. Here we only discuss fluid flows occupying the whole space so that the important effects of boundary layers are ignored. In most applications, $\nu$ is an extremely small quantity, typically ranging from $10^{-3}$ to $10^{-6}$ in turbulent flows. Thus one can anticipate that the behavior of inviscid solutions of the Euler equations with $\nu = 0$ is rather important in describing solutions of the Navier-Stokes equations when $\nu$ is small.

The local velocity gradient tensor solves

$$\partial_t M + (u \cdot \nabla)M + M^2 = \nu \Delta M - (\nabla \otimes \nabla)p. \tag{2.10}$$

Taking the trace of M and noting $trM = \nabla \cdot u = 0$ one has

$$trM^2 = -\Delta p. \tag{2.11}$$

This gives $p = -\Delta^{-1}(trM^2)$. The second term in (2.10) therefore amounts to the $n \times n$ time-dependent matrix
$$(\nabla \otimes \nabla)\Delta^{-1}(trM^2) = R[trM^2].$$
Here $R[w]$ denotes the $n \times n$ matrix whose entries are given by $(R[w])_{ij} := R_i R_j(w)$ where $R_j$ denote the Risez transforms, $R_j = -(-\Delta)^{-1/2}\partial_j$, i.e.,
$$\widehat{[R_j w]}(\xi) = -i\frac{\xi_j}{|\xi|}\hat{w}(\xi) \quad \text{for} \quad 1 \le j \le n.$$
This yields the equivalent formulation of NS equations which reads
$$(2.12) \qquad \partial_t M + (u \cdot \nabla)M + M^2 = \nu \Delta M + R[trM^2]$$
subject to the trace-free initial data
$$M(\cdot, 0) = M_0, \quad trM_0 = 0.$$
Note that the invariance of incompressibility is already taken into account in (2.12) since $\partial_t trM = 0$ and hence $trM = trM_0 = 0$. The inviscid case $\nu = 0$ in (2.12) gives the corresponding Euler equation. It is the *global* term in the above equations, $R[trM^2]$, which makes the problem rather intricate to solve, both analytically and numerically. Various simplifications to this pressure Hessian were sought, see, e.g. [34, 12, 5, 11].

Here we focus our attention on the so called restricted Euler equations proposed in [34] as a *localized* alternative to (2.12).

Specifically, we consider a gradient flow, $M$, governed by
$$(2.13) \qquad \partial_t M + (u \cdot \nabla)M + M^2 = \frac{1}{n}trM^2 I_{n \times n}.$$
We observe that as in the global model, the incompressibility is still maintained in this localized model, since $trM^2 = tr[trM^2 I_{n\times n}/n]$ implies $\partial_t trM = 0$.

For arbitrary $n \ge 3$, we use the spectral dynamics of $M$ in order to show a finite time breakdown of (2.13), generalizing the previous result of [34]. The finite time breakdown follows in §6 once we identify a set of $[n/2]+1$ global invariants in terms of the eigenvalues of $M$. Moreover, the precise topology of the flow at the breakdown time is studied in §6. Finally, in the Appendix we study the spectral dynamics of the general $n \times n$ problem from yet another perspective, extending the study of traces, $trM^k$, $k = 1, \cdots, n$ initiated in [34] for the special $n = 3$ case.

**2.4. Restricted Euler-Poisson equations.** The Euler-Poisson equations
$$(2.14) \qquad \rho_t + \nabla \cdot (\rho u) = 0, \quad x \in \mathbb{R}^n, \quad t \in \mathbb{R}^+,$$
$$(2.15) \qquad (\rho u)_t + \nabla \cdot (\rho u \otimes u) = k\rho \nabla \phi,$$
$$(2.16) \qquad \Delta \phi = \rho, \quad x \in \mathbb{R}^n,$$
are the usual statements of the conservation of mass, Newton's second law, and the Poisson equation defining, say, the electric field in terms of the charge. Here $k$ is a scaled physical constant, which signifies the property of the underlying forcing, repulsive if $k > 0$ and attractive if $k < 0$. The unknowns are the local density $\rho = \rho(x, t)$, the velocity field $u = u(x, t)$, and the potential $\phi = \phi(x, t)$.

If follows that
$$\partial_t u + u \cdot \nabla u = k\nabla \phi,$$

where the forcing $F = k\nabla\phi$ is the gradient of potential governed by the Poisson equation (2.16). Differentiation yields a local velocity gradient tensor which solves

$$\partial_t M + u \cdot \nabla M + M^2 = k(\nabla \otimes \nabla)\phi = kR[\rho],$$

where the coupling enters through the global term $kR[\rho]$, with density $\rho$ governed by

$$\partial_t \rho + u \cdot \nabla \rho + \rho tr M = 0.$$

Passing to Lagrangian coordinates, that is, using the change of variables $\alpha \mapsto x(\alpha, t)$ with $x(\alpha, t)$ solving

$$\frac{dx}{dt} = u(x, t), \quad x(\alpha, 0) = \alpha,$$

then Euler-Poisson equations become the coupled system

(2.17) $$\frac{d}{dt}M + M^2 = kR[\rho], \quad \frac{d}{dt} := \partial_t + u \cdot \nabla,$$

(2.18) $$\frac{d}{dt}\rho + \rho tr M = 0,$$

subject to initial condition

$$(M, \rho)(\cdot, 0) = (M_0, \rho_0).$$

Again, it is the nonlocal term, $R[\rho]$, which is the main obstacle, in the multi-dimensional setting $n > 1$, in the investigation of the Euler-Poisson system, see e.g. [21].

In this paper we introduce the corresponding restricted Euler-Poisson system

(2.19) $$\partial_t M + u \cdot \nabla M + M^2 = \frac{k}{n}\rho I_{n \times n},$$

(2.20) $$\partial_t \rho u \cdot \nabla \rho + \rho tr M = 0,$$

subject to initial data

$$(M, \rho)(\cdot, 0) = (M_0, \rho_0).$$

In §7 we use the spectral dynamics of $M$ in order to study the time regularity for this restricted Euler-Poisson model. Here we give a sufficient condition for the global existence of the 2D solutions which applies, for example, for a class of initial configurations with sufficiently large vorticity $|\omega_0| >> 1$. With other initial configurations, however, the finite time breakdown of solutions may – and actually does occur. Hence global regularity depends on whether the initial configuration crosses an intrinsic, $O(1)$ critical threshold, and we refer to [29] for a detailed study of the 2D critical threshold phenomena in this case. Moreover, for arbitrary $n \geq 3$ we obtain a family of $[n/2]$ global invariants, interesting for their own sake, with which the local topology of finite time breakdown is also characterized in §7.

## 3. Spectral dynamics of the Velocity Gradient Field

Let us rewrite the basic equation of velocity gradient field $M$ as

(3.1) $$\partial_t M + (u \cdot \nabla)M + M^2 = \nabla F,$$

where $\nabla F$ is a matrix involving spatial derivatives of the forcing.

It is usually difficult to quantify directly all entries in the velocity gradient tensor, $M$, and instead, suitable *linear* combinations like divergence, vorticity play a distinctive role in analysis. Here we show the special role played by the eigenvalues of the velocity gradient tensor, $\lambda(M)$, in governing the entries of $M$, and we note in passing, the strong nonlinear dependence of $\lambda(M)$ on the entries of $M$. Consult, for example, the nonintuitive $L^1$-contraction for 2D dusty medium model derived in (5.9) below.

The following lemma is in the heart of matter.

**Lemma 3.1** (Spectral dynamics). *Consider the general dynamical system (3.1) associated with arbitrary velocity field $u$ and forcing $F$. Let $\lambda(M)$ be a (possibly complex) eigenvalue of $M$ associated with corresponding left(right) eigenvector $l(r)$. Then the dynamics of $\lambda(M)$ is governed by the corresponding Ricatti-like equation*

$$\partial_t \lambda + u \cdot \nabla \lambda + \lambda^2 = <l, \nabla F r>.$$

*Remark* 3.2. If $F = 0$ one has the same equation for $\lambda$ as for $M$ with time-independent eigenvectors, thus $M(t)$ are isospectral. The difficulty lies in the eigenstructure induced by the forcing $<l, \nabla F r>$.

*Proof.* Let the left(right) eigenvectors of $M$ associated with $\lambda$ be $l(r)$, normalized so that $lr = 1$. Then one has

$$Mr = \lambda r, \quad lM = \lambda l.$$

Differentiation of the first relation with respect to $t$ gives

$$\partial_t M r + M \partial_t r = \partial_t \lambda r + \lambda \partial_t r.$$

Multiply $l$ on the left of the above equation to obtain

$$l \partial_t M r + \lambda l \partial_t r = \partial_t \lambda + \lambda l \partial_t r,$$

hence

$$l \partial_t M r = \partial_t \lambda.$$

Similarly differentiation of the relation $Mr = \lambda r$ with respect to $x_j$ leads to

$$\partial_j M r + M \partial_j r = \partial_j \lambda r + \lambda \partial_j r.$$

Multiply on the left by $lu_j$ with $lr = 1$ to get

$$lu_j \partial_j M r = u_j l \partial_j \lambda r = u_j \partial_j \lambda.$$

Therefore

$$lu \cdot \nabla M r = u \cdot \nabla \lambda.$$

A combination of the above facts together with $lM^2 r = \lambda^2$ gives

$$\partial_t \lambda + u \cdot \nabla \lambda + \lambda^2 = <l, \nabla F r>.$$

This completes the proof. □

## 4. Critical Thresholds for Linear Damping

Consider the convective equation

$$\partial_t u + u \cdot \nabla u = Cu, \quad u(x,0) = u_0(x),$$

with a simple linear damping $C$ being a constant matrix. The gradient tensor satisfies

$$\partial_t M + u \cdot \nabla M + M^2 = CM,$$

which suggests that the eigenvalues solve

(4.1) $$\partial_t \lambda + u \cdot \nabla \lambda + \lambda^2 = c\lambda,$$

where $c(t) = c_M(t) := lCr$. Here $l(r)$ are the left(right) eigenvectors of $M$ associated with the eigenvalue $\lambda$. Along the particle path $x = x(\alpha, t)$, defined by

$$\frac{d}{dt} x(\alpha, t) = u(t, x(\alpha, t)), \quad x(\alpha, 0) = \alpha, \quad \alpha \in \mathsf{IR}^n,$$

the Ricatti-type $\lambda$-equation amounts to

$$\frac{d}{dt} \lambda + \lambda^2 = c(t)\lambda.$$

The solution can be expressed in terms of $c(t)$ as

$$\lambda(t) = \frac{\lambda(0)b(t)}{1 + \lambda(0) \int_0^t b(\tau)d\tau}, \quad b(t) := \exp\left(\int_0^t c(\tau)d\tau\right).$$

¿From the above formula it follows that

**Lemma 4.1.** *Consider the eigenvalue equation (4.1) with initial data $\lambda(0)$.*
  *(1) If $Im(\lambda(0)b(t)) \neq 0$, then its solution remains regular for all time;*
  *(2) If $Im(\lambda(0)b(t)) = 0$, then its solution remains bounded as long as*

(4.2) $$Re\left(\lambda(0) \int_0^t b(\tau)d\tau\right) > -1.$$

For the simple example of a scalar damping, $C = -\beta I_{n \times n}$, $\beta > 0$, one has

$$\frac{d\lambda}{dt} + \lambda^2 = -\beta\lambda,$$

with a solution (corresponding to $b(t) = e^{-\beta t}$) given by

$$\lambda(t) = \frac{\lambda(\alpha, 0)e^{-\beta t}}{1 + \lambda(\alpha, 0)\beta^{-1}(1 - e^{-\beta t})}.$$

This solution is bounded from below for all time if and only if $\lambda(\alpha, 0)$ is either complex or

$$\inf_{\alpha \in \mathsf{IR}^n} \lambda(\alpha, 0) \geq -\beta,$$

which is a very simple form of a critical threshold phenomena. For more general examples of critical threshold phenomena, consult [14, 15, 29] and the study in §7 below.

## 5. Irrotational viscous flow

In this section we deal with viscous dusty medium flow, $u := u^\nu$, governed by

(5.1) $$\partial_t u + u \cdot \nabla u = \nu \Delta u, \quad u(x,0) = u_0(x).$$

The velocity gradient tensor $M := \nabla u$ satisfies

(5.2) $$\partial_t M + u \cdot \nabla M + M^2 = \nu \Delta M.$$

It follows that if the initial velocity is irrotational, $\nabla \times u_0 = 0$, then the flow remains irrotational, $\nabla \times u = 0$.

**Lemma 5.1** (Viscous Spectral Dynamics). *Assume that the flow is irrotational $\nabla \times u_0 = 0$. Then the real eigenvalues $\lambda = \lambda(\nabla u)$ satisfy*

$$\partial_t \lambda + u \cdot \nabla \lambda + \lambda^2 = \nu \Delta \lambda + Q.$$

*Here $Q$ satisfies the constraint*

$$a(\lambda_{min} - \lambda) \leq Q \leq a(\lambda_{max} - \lambda), \quad \lambda_{\{\max \atop \min\}} := \left\{ {\max \atop \min} \right\} \lambda(\nabla u),$$

*where $a$ is given by*

$$a := 2\nu \sum_k \partial_k r^\top \partial_k r > 0$$

*and $r$ being the right eigenvector of $\nabla u$ associated with $\lambda$.*

*Proof.* Let $l(r)$ be the normalized left(right) eigenvectors of $M$ associated with the eigenvalue $\lambda$, then one has

$$\partial_t \lambda + u \cdot \nabla \lambda + \lambda^2 = \nu l \Delta M r.$$

Observe that $M$ is symmetric due to the fact that $\nabla \times u = 0$, and consequently $\lambda$ are all real quantities. Differentiation of $lM = \lambda l$ with respect to $x$ twice gives

$$\Delta l M + 2 \nabla l \cdot \nabla M + l \Delta M = \Delta \lambda l + 2 \nabla \lambda \cdot \nabla l + \lambda \Delta l,$$

which upon multiplication against $r$ on the right leads to

$$l \Delta M r = \Delta \lambda + 2 \left[ (\nabla \lambda \cdot \nabla l) r - (\nabla l \cdot \nabla M) r \right].$$

Here the differentiation operators apply component wise, e.g., $\nabla l \cdot \nabla M = \sum_{k=1}^n \partial_k l \partial_k M$. On the other hand it follows from $Mr = \lambda r$ that

$$\nabla M r = \nabla \lambda r + \lambda \nabla r - M \nabla r.$$

This gives

$$(\nabla l \cdot \nabla M) r = \nabla l \cdot \nabla \lambda r + \lambda \nabla l \cdot \nabla r - \nabla l \cdot M \nabla r.$$

A combination of the above facts yields

$$Q = 2\nu \left[ -\lambda \sum_{k=1}^n \partial_k l \partial_k r + \sum_{k=1}^n \partial_k l M \partial_k r \right].$$

Since the flow is irrotational we have $M^\top = M$ and $l = r^\top$, with upper-index $^\top$ denoting the transpose. The second term in $Q$ is then bounded by

$$\lambda_{min} \sum_{k=1}^n \partial_k r^\top \partial_k r \leq \sum_{k=1}^n \partial_k l M \partial_k r \leq \lambda_{max} \sum_{k=1}^n \partial_k r^\top \partial_k r,$$

which completes the proof.  □

Here the question of interest for us is the convergence of $u = u^\nu$ as $\nu \to 0$. To answer this question, it suffices to show the precompactness of the family $\{u^\nu\}_{\nu>0}$. It is here that we take advantage of the spectral dynamics of the velocity gradient tensor, $\nabla u$. For the 2D case we shall show the precompactness via several lemmata. We start with the essential

**Lemma 5.2** ($L^1$-Contraction). *Let $\lambda_i$, $i = 1, 2$, be two (real) eigenvalues of velocity gradient field $\nabla u^\nu$ in (5.1). If $(\lambda_2 - \lambda_1)(0) \in L^1(\mathbb{R}^2)$, then*
$$\|(\lambda_2 - \lambda_1)(t)\|_{L^1(\mathbb{R}^2)} \leq \|(\lambda_2 - \lambda_1)(0)\|_{L^1(\mathbb{R}^2)}.$$

*Proof.* In the 2D case we have $\lambda_{min} = \lambda_1 \leq \lambda_2 = \lambda_{max}$. Setting $\eta = \lambda_2 - \lambda_1$ one has from Lemma 4.1
$$\partial_t \eta + u \cdot \nabla \eta + \eta(\lambda_1 + \lambda_2) \leq \nu \Delta \eta.$$
Observe that $\nabla \cdot u = \lambda_1 + \lambda_2$ which yields
$$\partial_t \eta + \nabla \cdot (\eta u) \leq \nu \Delta \eta.$$
Spatial integration gives the $L^1$ estimate for $\eta = \lambda_2 - \lambda_1 \geq 0$ as asserted.  □

Next, from Lemma 5.1 we see that the largest eigenvalue $\lambda_{\max}$ satisfies the differential inequality
$$\partial_t \lambda_{\max} + u \cdot \nabla \lambda_{\max} + \lambda_{\max}^2 \leq \nu \Delta \lambda_{\max},$$
and by a comparison principle we obtain
$$(5.3) \qquad \lambda_{\max}(t) \leq \frac{1}{\lambda(0)^{-1} + t} \leq \frac{1}{t}.$$

We note in passing that this, combined with the symmetry of gradient field $\nabla u$, is equivalent to the one-sided entropy-type estimate $\sup_{\|\xi\|=1} \xi^\top \nabla u \xi \leq 1/t$, which coincides with the well known semi concavity property in the context of convex Hamilton-Jacobi equations, see, e.g., [25, 28, 32].

The above one-sided bounds enable us to establish the following.

**Lemma 5.3** (*BV* Bound). *Consider the dusty medium equation (5.1) with compactly supported initial data $u_0^\nu = u^\nu(x, 0)$ such that $\|u_0^\nu\|_{BV(\mathbb{R}^2)}$ is bounded uniformly in $\nu$. Then the corresponding velocity, $u^\nu$, satisfies*
$$\|u^\nu(\cdot, t)\|_{BV(\mathbb{R}^2)} \leq \text{Const}.$$
*Moreover, for $t_1, t_2 \geq 0$ we also have*
$$(5.4) \qquad \|u^\nu(x, t_2) - u^\nu(x, t_1)\|_{L^1(\mathbb{R}^2)} \leq \text{Const}.|t_2 - t_1|^{1/3}.$$

*Proof.* The one-sided upper bound for $\lambda_{max}$, (5.3), implies that the positive part of the divergence, $(u_x + v_y)_+ = (\lambda_1 + \lambda_2)_+$ is bounded. We observe that $\lambda_1, \lambda_2$ are essentially supported on a finite domain in the sense of their exponential decay outside a finite region of propagation, and hence $\int_{\mathbb{R}^2} (u_x + v_y)_+ \leq \text{Const}$. This, combined with $\int_{\mathbb{R}^2} u_x + v_y = 0$, yields that $\lambda_1 + \lambda_2 = u_x + v_y \in L^1(\mathbb{R}^2)$. Augmented with the fact that $\lambda_2 - \lambda_1 \in L^1(\mathbb{R}^2)$ we conclude
$$(5.5) \qquad \lambda_i \in L^1(\mathbb{R}^2), \quad i = 1, 2.$$

This gives
$$\int_{\mathbb{R}^2} \|\nabla u^\nu\| dxdy = \int_{\mathbb{R}^2} \|diag(\lambda_1, \lambda_2)\| dxdy < \infty,$$
with the usual matrix norm, $\|\cdot\|$, defined as $\|M\| =: \sup_{\|\xi\|=1} |M\xi|$. In fact, since $M = \nabla u$ is symmetric, there exists a unitary matrix $U$ such that $U^\top M U = diag(\lambda_1, \lambda_2)$, and hence
$$\|\nabla u\| = \|U^\top M U\| = \|diag(\lambda_1, \lambda_2)\|.$$
Thus, the $BV$ bound of $u^\nu$ follows from (5.5). To estimate the modulus of continuity in time, we multiply equation (5.1) by a smooth test function $\psi \in C_0^\infty$ and use the spatial BV bound to obtain
$$\left| \int_{\mathbb{R}^2} \psi(x)(u(x,t_2) - u(x,t_1))dx \right| \leq \text{Const.}(t_2-t_1)(|\psi|_\infty + |\Delta\psi|).$$
This inequality and the BV estimate combined with Kružkov's interpolation theorem [23, page 233] yield (5.4). □

In order to identify the vanishing viscosity limit, $\lim_{\nu \to 0} u^\nu$, we introduce the notion of a weak solution for corresponding inviscid equation

(5.6) $$\partial_t u + u \cdot \nabla u = 0.$$

For irrotational flow, $\nabla \times u = 0$, one has $u \cdot \nabla u \equiv \nabla(|u|^2/2)$, and the reduced inviscid equation (5.6) can be recast into the conservative form
$$\partial_t u + \nabla \left(\frac{|u|^2}{2}\right) = 0.$$

The irrotational property of both viscous and inviscid flows suggests that there exists a potential $\phi$ such that $u = \nabla \phi$, where $\phi$ solves the Hamilton-Jacobi equation

(5.7) $$\phi_t + \frac{1}{2}|\nabla \phi|^2 = 0, \quad \phi(x,0) = \phi_0.$$

According to the classical theory of the Hamilton-Jacobi equation [7, 8], there exists a unique continuous solution $\phi(x,t)$ to the above problem, expressed via the Hopf-Lax's formula, [13, page 560], $\phi(x,t) = \min_{y \in \mathbb{R}^n} \{t|x-y|^2/2 + \phi_0(y)\}$. We make

**Definition 5.1.** A measurable function $u$ is called a weak solution of the inviscid equation (5.6) if $u = \nabla \phi$ with the potential $\phi$ being the unique weak solution of the Eikonal equation (5.7).

Equipped with this definition of a weak solution, we now turn to summarize our convergence results by stating

**Theorem 5.4** (Vanishing viscosity limit)**.** *Consider the dusty medium equation (5.1) with irrotational initial data $u^\nu(\cdot,0) \in L^1 \cap L^\infty(\mathbb{R}^2)$ such that*
$$u^\nu(x,0) \to u_0(x) \quad in \quad L^1(\mathbb{R}^2).$$
*Then, the local velocity $u^\nu$ converges to the unique weak solution of (5.6), i.e., we have*

(5.8) $$u^\nu(x,t) \to u(x,t) \quad in \quad L^\infty([0,T]; L^1(\mathbb{R}^2)),$$

*where $u = \nabla \phi$ is the viscosity solution of the Eikonal equation (5.6).*

*Proof.* We begin by first assuming that $u^\nu(x,0)$ is compactly supported in $BV(\mathbb{R}^2)$, uniformly with respect to $\nu$. By Lemma 4.2, $u^\nu$ have uniformly bounded spatial variation, i.e.,
$$\|u^\nu(\cdot,t)\|_{BV(\mathbb{R}^2)} \leq \text{Const}.$$
Hence $\{u^\nu(x,t), 0 \leq t \leq T\}$ is a bounded set in $L^1 \cap BV(\mathbb{R}^2)$ and by Helly's theorem it is therefore precompact in $L^1_{loc}(\mathbb{R}^2)$. Note that $\|u^\nu(x,t)\|_{L^1(\mathbb{R}^2)}$ is Hölder continuous in time, and by Cantor diagonalization process of passing to further subsequence if necessary, (5.8) follows. This completes the corresponding proof for compactly supported $BV$ initial data. The general case is justified by standard cutoff and BV-regularization of arbitrary $L^1 \cap L^\infty(\mathbb{R}^2)$ initial data.

It remains to show that the limit function $u(x,t)$ satisfy the weak formulation. It follows from the equation for $u^\nu$ that
$$\partial_t u^\nu + \nabla\left(\frac{1}{2}|u^\nu|^2\right) = \nu \Delta u^\nu.$$
We multiply this identity by $\psi(x,t) \in C_0^\infty(\mathbb{R}^2)$ and integrate by parts to get
$$\int_{R^2}\left[-\psi_t u^\nu - \nabla\psi \frac{|u^\nu|^2}{2}\right] dxdy = \nu \int_{\mathbb{R}^2} u^\nu \Delta\psi\, dxdy.$$
Note that $W^{1,1}$ is embedded into $L^2(\mathbb{R}^2)$ for the two-dimensional case. Thus passing to the limit $\nu \to 0$ one obtains the desired weak formulation. $\square$

*Remark* 5.5. We would like to point out that the above convergence result can be obtained at the level of Hamilton-Jacobi equations. The equivalence between the weak entropy solutions to conservation laws and the viscous solutions to the corresponding Hamilton-Jacobi equations has been known in the literature, see e.g. [4, 7, 24, 25, 22]. The point made here is that we obtain the compactness at the level of $u$ by using the spectral dynamics of its velocity gradient tensor, $\nabla u$, which is independent of a maximum principle at the level of HJ equations. In particular, the $2D$ $L^1$-contraction stated in Lemma 5.2, recast at the level of HJ equation (5.7), amounts to the nonintuitive apriori estimate

(5.9) $$\left\|\sqrt{(\Delta\phi)^2 - J(\phi)}(\cdot,t)\right\|_{L^1} \leq \left\|\sqrt{(\Delta\phi_0)^2 - J(\phi_0)}\right\|_{L^1}, \quad J(\phi) := \phi_{xx}\phi_{yy} - \phi_{xy}^2.$$

## 6. Restricted Euler Dynamics

### 6.1. Spectral dynamics and global invariants.
We now turn to discuss the restricted Euler dynamics, which is a localized version of the full Euler/Navier-Stokes equation (2.12):
$$\partial_t M + (u \cdot \nabla)M + M^2 = \nu \Delta M + R[tr M^2].$$
By the definition of the operator $R$, one has
$$R[tr M^2] = \nabla \otimes \nabla \Delta^{-1}[tr M^2] = \nabla \otimes \nabla \int_{\mathbb{R}^n} K(x-y)(tr M^2)(y)dy,$$
where the kernel $K(\cdot)$ is given by
$$K(x) = \begin{cases} \frac{1}{2\pi}\ln|x| & n=2, \\ \frac{1}{(2-n)\omega_n|x|^{n-2}} & n>2, \end{cases}$$

with $\omega_n$ denoting the surface area of the unit sphere in $n$-dimensions. A direct computation yields

$$\partial_i \partial_j K * trM^2 = \frac{trM^2}{n}\delta_{ij} + \int_{\mathbb{R}^n} \frac{|x-y|^2 \delta_{ij} - n(x_i - y_i)(x_j - y_j)}{\omega_n |x-y|^{n+2}} trM^2(y) dy.$$

This shows that the local part of the global term $R[trM^2]$ is $trM^2/nI_{n\times n}$. We thus use this local term, $trM^2/nI_{n\times n}$, to approximate the pressure Hessian. The corresponding local gradient field then evolves according to the following restricted Euler model

(6.1) $$\partial_t M + u \cdot \nabla M + M^2 = trM^2/nI_{n\times n}.$$

This is a matrix Ricatti equation which, as we shall see below, is responsible for the formation of singularities at finite time. We note that with this local model, all particles evolve independently of each other. The mixing due to the global forcing in the general Euler dynamics, however, could prevent this type of finite time breakdown.

Nevertheless, as a local approximation of the pressure Hessian, the above model, so called restricted Euler dynamics, has caught great attention since first introduced in [34], because it can be used to understand the local topology of the Euler dynamics and to capture certain statistical features of the physical flow.

Consider a bounded, divergence-free, smooth vector field $u : \mathbb{R}^n \times [0, T] \to \mathbb{R}^n$. Let $x = x(\alpha, t)$ denote an orbit associated to the flow by

$$\frac{dx}{dt} = u(x, t), \quad 0 < t < T, \quad x(\alpha, 0) = \alpha \in \mathbb{R}^n.$$

Then along this orbit, the velocity gradient tensor of the restricted Euler equations (6.1) satisfies

$$M' + M^2 = trM^2/nI_{n\times n}, \quad ' := \frac{d}{dt}.$$

By the spectral dynamics lemma 3.1, the corresponding eigenvalues satisfy

(6.2) $$\lambda_i' + \lambda_i^2 = \sum_{k=1}^{n} \lambda_k^2/n, \quad i = 1, \cdots, n.$$

This is a closed system which serves as a simple approximation for the evolution of the velocity gradient field.

Let us start by revisiting the case $n = 3$, consult [34, 5], for which we will present below a complete phase-plane analysis expressed in terms of $\lambda_i's$. Subtraction of two consecutive equations in (6.2) gives the following equivalent system

$$[ln(\lambda_1 - \lambda_2)]' + \lambda_1 + \lambda_2 = 0,$$
$$[ln(\lambda_2 - \lambda_3)]' + \lambda_2 + \lambda_3 = 0,$$
$$[ln(\lambda_3 - \lambda_1)]' + \lambda_3 + \lambda_1 = 0.$$

Summation of these three equations and taking into account the incompressibility condition, $\sum_{i=1}^{3} \lambda_i = 0$, yields the following global invariant

$$(\lambda_1 - \lambda_2)(\lambda_2 - \lambda_3)(\lambda_3 - \lambda_1) = \text{Const.}$$

This invariant projected onto the phase plane $(\lambda_1, \lambda_2)$, recast into

$$(\lambda_2 - \lambda_1)(\lambda_2 + 2\lambda_1)(2\lambda_2 + \lambda_1) = \text{Const.},$$

which serves as a global invariant of the 2 × 2 system

(6.3) $$\lambda_1' = [-\lambda_1^2 + 2\lambda_2^2 + 2\lambda_1\lambda_2]/3,$$
(6.4) $$\lambda_2' = [2\lambda_1^2 - \lambda_2^2 + 2\lambda_1\lambda_2]/3.$$

We then have three separatrixes passing through the origin, which is the only rest point in this case,

$$\lambda_1 = \lambda_2, \quad \lambda_1 = -2\lambda_2 \leftrightarrow (\lambda_2 = \lambda_3), \quad \lambda_1 = -\frac{1}{2}\lambda_2 \leftrightarrow (\lambda_3 = \lambda_1).$$

The vector field in the phase plane is drawn in Figure 6.1.

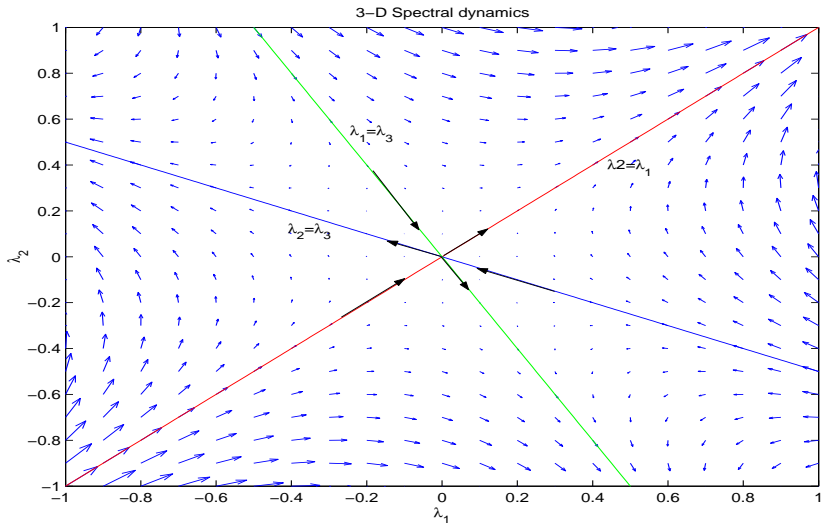

FIGURE 6.1. 3-D Spectral dynamics of the Restricted Euler Equations

Three special solutions corresponding to the separatrixes can be obtained explicitly. Consider, for example, the separatrix $\lambda_1 = \lambda_2$, for which $\lambda_1$ is necessarily a real solution of the Ricatti equation

$$\lambda_1' = \lambda_1^2.$$

The solution, given by

(6.5) $$\lambda_1(x,t) = \frac{\lambda_1(\alpha,0)}{1 - \lambda_1(\alpha,0)t},$$

is bounded if and only if the real $\lambda_1$ is nonpositive, $\lambda_1(\alpha, 0) \leq 0$.

Next if $\lambda_i(0)$ is complex then $\lambda_i$ remains complex later on. Let $(\lambda_1, \lambda_2)$ be a complex pair with $\lambda_2 = \bar{\lambda}_1 = \gamma + \beta i$, then (6.3), (6.4) recast into

(6.6) $$\beta' = -2\beta\gamma, \quad \gamma' = \gamma^2 + \frac{1}{3}\beta^2.$$

Solving the above $2 \times 2$ system gives the following invariant, $(\beta^2 + 9\gamma^2)\beta = $ Const; it follows that the general solution passing the rest point $(0,0)$ must be real,

$$\beta = 0, \quad \gamma(t) = \frac{\gamma(0)}{1 - \gamma(0)t},$$

which is reduced back to the first case in (6.5). Note that if the eigenvalues are complex, the Lagrangian trajectories are rotating, and if the eigenvalues are all real, the Lagrangian trajectories are just strain dominated, see Figure 6.2.

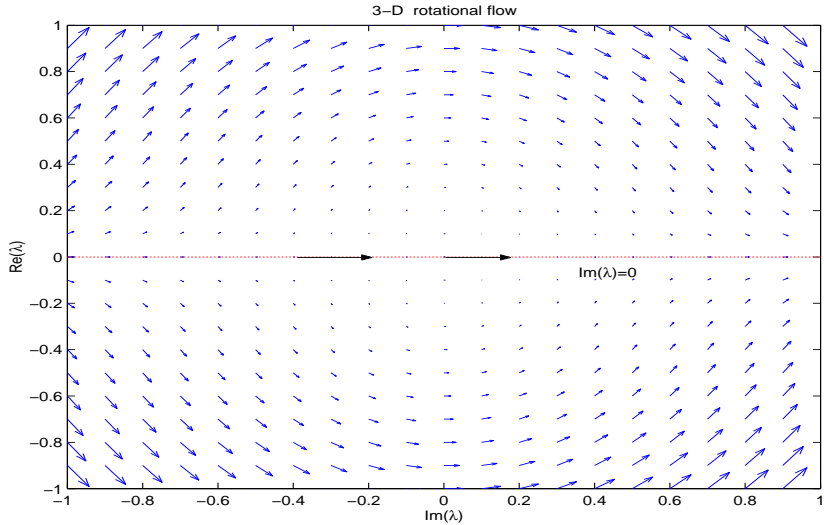

FIGURE 6.2. 3-D Rotational flow in Restricted Euler Equations

We now summarize by stating the following

**Lemma 6.1** (Topology of flow in 3-D restricted Euler). *A global invariant of the 3D restricted Euler equations (6.2) is given by*

(6.7) $$(\lambda_1 - \lambda_2)(\lambda_2 - \lambda_3)(\lambda_3 - \lambda_1) = \text{Const}.$$

*The three explicit solutions passing through origin are*

$$(\lambda_1, \lambda_2, \lambda_3)(t) = \{(1,1,-2), (1,-2,1), (-2,1,1)\} \frac{a(\alpha,0)}{1 - a(\alpha,0)t}$$

*All other solutions will develop finite time singularity. If an eigenvalue is complex, then the Lagrangian trajectories are rotating.*

Lemma 6.1 deals with the 3-dimensional restricted Euler equations which were studied earlier in [34] using a different approach based on trace dynamics, consult the Appendix below. Here we enjoy the advantage of being able to generalize our spectral dynamics approach taken in Lemma 6.1 to the arbitrary $n$-dimensional case. The global invariants in such $n$-dimensional systems are tied to a particular set of sequences of indices $\mathcal{I} = $ $\{\mathcal{I}_1, \mathcal{I}_2, \ldots\}$ with each $\mathcal{I}_k$ being a sequence of pairs of different indices, $(i \neq j)$, such that there exists an integer $N := N(n)$ for which

$$\text{(6.8)} \qquad \sum_{(i,j) \in \mathcal{I}} (\lambda_i + \lambda_j) = N \sum_{k=1}^n \lambda_k, \quad \forall \lambda's.$$

There are several ways of forming these pairs, $(i,j)$, so that (6.8) holds. Here is one:
• For even $n = 2m$ we let $\mathcal{I}_\sigma = (i,j) = \{(\sigma(2k-1), \sigma(2k))\}_{k=1,2,\ldots,m}$, ranging overall permutations $\sigma(\cdot)$ so that (6.8) holds with $N(n)_{|\{n \text{ even}\}} = 1$;
• For odd $n = 2m+1$ we let $\mathcal{I}_{\sigma\mu} = (i,j) = \{(\sigma(k) \neq \mu(k))\}_{k=1,2,\ldots,n}$, ranging overall permutations $\sigma(\cdot), \mu(\cdot)$ so that (6.8) holds with $N(n)_{|\{n \text{ odd}\}} = 2$.

The following lemma reveals the role such $\mathcal{I}$-pairs play in forming global invariants for the restricted Euler system (6.1).

**Lemma 6.2.** *[Global Invariants] Consider the n-dimensional restricted Euler system (6.1) subject to incompressible initial data, $\sum_{i=1}^n \lambda_i(0) = 0$. Then it admits the following global invariants in time:*

$$\text{(6.9)} \qquad \sum_{i=1}^n \lambda_i(t) = 0,$$

$$\text{(6.10)} \qquad \Pi_{(i,j) \in \mathcal{I}} (\lambda_i(t) - \lambda_j(t)) = \text{Const.}$$

*Proof.* Summation of the equations in (6.2) over index $i$ gives $[\sum_{i=1}^n \lambda_i(t)]' = 0$, which is combined with the incompressibility assumption $\sum_{i=1}^n \lambda_i(0) = 0$ to yield $\sum_{i=1}^n \lambda_i(t) = 0$.

For (6.10) we follow our previous argument in the 3D case. Subtracting the $j$-th equation from $i$-th equation in (6.2) yields

$$\text{(6.11)} \qquad [\alpha_{i,j}]_t + (\lambda_i + \lambda_j)\alpha_{i,j} = 0, \quad \alpha_{i,j} = \lambda_i - \lambda_j.$$

Divided by $\alpha_{i,j}$ and sum those equations in (6.11) with indices $(i,j) \in \mathcal{I}$, we have

$$\text{(6.12)} \qquad [ln \, \Pi \, \alpha_{i,j}]_t + N \sum_{k=1}^n \lambda_k = 0.$$

By incompressibility $\sum_{k=1}^n \lambda_k = 0$ and the global invariants asserted in (6.10) follow. $\square$

Two prototype examples are in order. In the 3D case we recover the global invariant (6.7), $(\lambda_1 - \lambda_2)(\lambda_2 - \lambda_3)(\lambda_3 - \lambda_1) = \text{Const.}$, corresponding to the sequence of pairs $\mathcal{I}_1 = \{(1,2), (2,3), (3,1)\}$. In the 4D case we have, in addition to incompressibility, the two global invariants, $\Pi_1 := (\lambda_1 - \lambda_2)(\lambda_3 - \lambda_4) = \text{Const}_1$ and $\Pi_2 := (\lambda_1 - \lambda_3)(\lambda_2 - \lambda_4) = \text{Const}_2$, corresponding to the two $\mathcal{I}$-sequences of indices, $\mathcal{I}_1 = \{(1,2), (3,4)\}$ and $\mathcal{I}_2 = \{(1,3), (2,4)\}$. Observe that a third global invariant corresponding to $\{(1,4), (1,3)\}$, is in fact generated by the difference of the first two, namely $\Pi_3 = (\lambda_1 - \lambda_4)(\lambda_2 - \lambda_3) \equiv \Pi_1 - \Pi_2$. Our next issue is therefore, a proper counting of these global invariants.

6.2. **On the number of global invariants.** The proof of Lemma 6.2 makes clear the direct linkage between each $\mathcal{I}$-sequence of indices, $(i,j)$ satisfying (6.8), and a global invariant formed by the corresponding product, $\Pi_{(i,j)}(\lambda_i - \lambda_j)$. Of course, not all the different $\mathcal{I}$-sequences satisfying (6.8) should be counted, since some of them lead to the same invariant products. We also need to remove any redundancy due to linear and nonlinear dependence among these different invariant products. Thus, we inquire about the following

*Question.* How many *independent* products, $\Pi_{(i,j)\in\mathcal{I}}(\lambda_i - \lambda_j)$ can be formed by $\mathcal{I}$-sequences, i.e., sequences of indices $(i,j)$ satisfying (6.8),

$$\exists N = N(n) \in \mathbb{Z} \text{ s.t. } \sum_{(i,j)\in\mathcal{I}}(\lambda_i + \lambda_j) = N\sum_{k=1}^{n}\lambda_k, \qquad \forall \lambda's?$$

We know that the number of such independent invariant products together with the incompressibility constraint (6.9) does not exceed $n$, the number of independent global invariants of the restricted Euler (6.1), and hence there are no more than $n-1$ independent invariants of form (6.10). But the precise answer remains open, and in particular we are not clear whether *all* the global invariants of (6.1) are necessarily the products formed in Lemma 6.2. Below we provide a *lower bound* for our question, by the construction of $\left[\frac{n}{2}\right]$ such independent invariants.

Let us begin by referring to the 4D example mentioned above. Starting with the first invariant, $\Pi_1 = (\lambda_1 - \lambda_2)(\lambda_3 - \lambda_4)$, we derive a second independent invariant by exchanging, $2 \leftrightarrow 3$, which leads to $\Pi_2 = (\lambda_1 - \lambda_3)(\lambda_2 - \lambda_4)$. Other possible exchanges are redundant, say $2 \leftrightarrow 4$ yields the linearly dependent product $\Pi_3 = (\lambda_1 - \lambda_4)(\lambda_2 - \lambda_3) \equiv \Pi_1 - \Pi_2$, and this is consistent with the fact that $i = 3, j = 4$ playing a symmetric role in the original $\Pi_1$-pair $(\lambda_3 - \lambda_4)$. We conclude that while forming the linearly independent products, $\Pi_{(i,j)\in\mathcal{I}}(\lambda_i - \lambda_j)$, at most one 'admissible' exchange between different pairs of $\mathcal{I}$-indices is permitted. Moreover, we should also exclude nonlinear dependence. For $n = 8$, for example, consider the four products, $\Pi_1 = (\lambda_1 - \lambda_2)(\lambda_3 - \lambda_4)(\lambda_5 - \lambda_6)(\lambda_7 - \lambda_8)$, $\Pi_2 = (\lambda_1 - \lambda_3)(\lambda_2 - \lambda_4)(\lambda_5 - \lambda_6)(\lambda_7 - \lambda_8)$, $\Pi_3 = (\lambda_1 - \lambda_2)(\lambda_3 - \lambda_4)(\lambda_5 - \lambda_7)(\lambda_6 - \lambda_8)$ and $\Pi_4 = (\lambda_1 - \lambda_3)(\lambda_2 - \lambda_4)(\lambda_5 - \lambda_7)(\lambda_6 - \lambda_8)$. The four invariant products are linearly independent — indeed, $\sum \alpha_k \Pi_k = 0$ with $\lambda_7 = \lambda_8$ is reduced to a linear combination of the last two 4D independent pairs, $\alpha_3 \Pi_3 + \alpha_4 \Pi_4 = 0 \implies \alpha_3 = \alpha_4 = 0$, and similarly, setting $\lambda_6 = \lambda_8$ yields $\alpha_1 = \alpha_2 = 0$. Nevertheless, they are redundant in view of their nonlinear dependence, $\Pi_4 \equiv \Pi_2 \times \Pi_3 / \Pi_1$.

Our construction of independent invariants in the general $n$-dimensional case proceeds as follows. We start, for even $n = 2m$, with the usual ordering $\mathcal{I}_1 = (1,2)(3,4)\ldots(n-1,n)$. Making an admissible exchange between the first and second pairs yields the next independent invariant associated with $\mathcal{I}_2 = (1,3)(2,4),(5,6)\ldots(n-1,n)$. Next, we make an admissible exchange between the second and third pairs, $\mathcal{I}_3 = (1,2)(3,5)(4,6)\ldots(n-1,n)$, and so on. In this manner we proceed with one admissible exchange between each two *consecutive* pairs, leading to the $m$ global invariants of the restricted Euler equations (9.1),

$$(6.13) \quad \begin{cases} \Pi_1 & := (\lambda_1 - \lambda_2)(\lambda_3 - \lambda_4) \cdot \ldots \cdot (\lambda_{n-1} - \lambda_n), \\ \Pi_k & := (\lambda_1 - \lambda_2) \cdot \ldots \cdot (\lambda_{2k-3} - \lambda_{2k-1})(\lambda_{2k-2} - \lambda_{2k}) \cdot \ldots \cdot (\lambda_{n-1}, \lambda_n), \\ & \quad k = 2, 3, \ldots, m. \end{cases}$$

To verify that these $m = \frac{n}{2}$ global invariants are independent, we note that by setting $\lambda_{2k-3} = \lambda_{2k-2}$ we have $\Pi_j = \begin{cases} \equiv 0, & j \neq k \\ \neq 0 & j = k \end{cases}$, which excludes the possible dependence $\Pi_k \neq \mathcal{F}(\Pi_1, \ldots, \Pi_{k\pm1}, \ldots, \Pi_n)$, $k = 1, 2, \ldots m$.

A similar procedure applies to the odd case, $n = 2m+1$. Starting with the usual ordering $\mathcal{I}_1 = (1,2), (2,3), \ldots (n,1)$, we make an admissible exchange between the first and third pairs, $\mathcal{I}_2 = (1,3)(2,3)(2,4)\ldots$, the third and fifth pairs, $(1,2)(2,3)(3,5)(4,5)(4,6)\ldots$ and so on, leading to the $m$ independent global invariants

$$(6.14) \quad \begin{cases} \Pi_1 &:= (\lambda_1 - \lambda_2)(\lambda_2 - \lambda_3)\cdot\ldots\cdot(\lambda_n - \lambda_1), \\ \Pi_k &:= (\lambda_1 - \lambda_2)\ldots(\lambda_{2k-3} - \lambda_{2k-1})(\lambda_{2k-2} - \lambda_{2k-1})(\lambda_{2k-2} - \lambda_{2k})\ldots(\lambda_n - \lambda_1), \\ & k = 2, 3, \ldots, m. \end{cases}$$

We conclude with

**Lemma 6.3.** *[Global Invariants] The $n$-dimensional restricted Euler system (6.1) subject to the incompressible initial data, $\sum_{i=1}^n \lambda_i(0) = 0$, admits the following $\left[\frac{n}{2}\right] + 1$ global invariants in time: the incompressibility (6.10), $\sum_{i=1}^n \lambda_i(t) = 0$, and the additional $n/2$ (– respectively, $(n-1)/2$) invariants specified in (6.13) for $n$ even (and respectively in (6.14) for $n$ odd).*

**6.3. Behavior at the finite breakdown time.** The rest of section is devoted to study the topology of the flow at the breakdown time based on the Lemma 6.2. We start by noting that the level set of the integrals of the restricted flow (6.9), (6.10) are not compact, and hence we have to perform singularity analysis to figure out in which orthant the flow may diverge. The idea is to build local solutions around the singularities in order to study the blow up-rate and the location where the finite-time blow-up actually occurs. The singularity analysis is a standard method to prove the integrability of ODEs. For readers' convenience we sketch the main steps below, and refer to [17] and references therein for more details of this method.

We assume a flow governed by the nonlinear ODE $w' = f(w)$, diverges at a finite time $t^*$, and we then seek local solutions of the form

$$w = \omega\tau^p\left[1 + \sum_{j=1}^\infty a_j \tau^{j/q}\right],$$

where $\tau = t^* - t$, $p \in \mathbb{R}^n$, $q \in \mathbb{N}$ and $a_j$ is a polynomial in $\log(t^* - t)$ of degree $N_j \leq j$. There are three steps to determine the above series : (1) find the so called balance pair, $(\omega, p)$, such that the dominant behavior, $\omega\tau^p$, is an exact solution of some truncated system $w' = \tilde{f}(w)$; (2) computation of the resonances, which are given by the eigenvalues of the matrix $-\frac{\partial \tilde{f}(w)}{\partial w} - diag(p)$; (3) the last step of the singularity analysis consists of finding the explicit form for the different coefficients $a_j$ by inserting the full series in the original system, $w' = f(w)$.

Armed with the above algorithm, we proceed to carry out the singularity analysis for the restricted Euler equations. Let the dominant behavior of the $\lambda$-system (6.2) assume the form

$$\lambda_i(t) \sim \omega_i \tau^{p_i}, \quad i = 1, \cdots, n.$$

Upon substitution into (6.2) one has

$$-\omega_i p_i \tau^{p_i - 1} + \omega_i^2 \tau^{2p_i} = \frac{1}{n}\sum_{j=1}^n \omega_j^2 \tau^{2p_j}.$$

Equating the powers of $\tau$ as $\tau \to 0$ we find, $p_i = -1$, and the $\omega_i$'s satisfy the equation

$$\omega_i + \omega_i^2 = \frac{1}{n}\sum_{j=1}^{n}\omega_j^2.$$

There is a n-parameter family of such $\omega$'s

$$\omega^{(k)} = \left(\frac{1}{n-2}, \cdots, \frac{1-n}{n-2}, \cdots, \frac{1}{n-2}\right), \quad k = 1, \cdots, n.$$

Due to the symmetry of the equation, the flow may diverge in $n$ out of $2^n$ orthants. More precisely, we have

**Lemma 6.4.** *The only $n$ stable solutions of the spectral dynamics (6.2), $\Lambda = (\lambda_1, \cdots, \lambda_n)$, associated with restricted Euler equations (6.1) are explicitly given by*

$$\Lambda = \Lambda^{(k)}(x,t) = \omega^{(k)}\frac{(n-2)a(x)}{n-2-a(x)t}, \quad \Lambda^{(k)}(0) = \omega^{(k)}a(x), \quad k = 1, \cdots, n,$$

*with arbitrary $a(x) \leq 0$.*

To sum up, we state the following

**Theorem 6.5.** *Consider the restricted Euler dynamics (6.2) with initial data $(\lambda_1(0), \cdots, \lambda_n(0))$. The level set of $\left[\frac{n}{2}\right]$ global invariants given by*

$$\amalg_{i,j\in\mathcal{I}}(\lambda_i - \lambda_j) = \text{Const}$$

*are not compact. The general solution may break down at finite time in one of the $n$ orthants $\{+, +, \cdots, -, \cdots +, +\}$ along the $k$th separatrix*

$$(1, \cdots, 1-n, \cdots, 1)\frac{a(x)}{n-2-a(x)t}, \quad k = 1, \cdots, n$$

*whenever $a(x) > 0$.*

*Remark* 6.6. Other possible variants of the local restricted Euler equations can be written in the form

$$\frac{d}{dt}M + \theta(M^2 - trM^2/nI_{n\times n}) = 0$$

with $\theta \in (0, \infty)$. This equation becomes anisotropic, but the local topology of the solution remains the same as in the isotropic model (6.1) below. Indeed a hyperbolic scaling, $(t, x) \to (\theta t, \theta x)$, leads to the isotropic model corresponding to $\theta = 1$.

## 7. Restricted Euler-Poisson dynamics

We begin by introducing the so called restricted Euler-Poisson equations. As argued in §6 we retain the local part of the nonlocal term $kR[\rho]$ in the Lagrangian form of the Euler-Poisson (EP) system (2.17), (2.18) to obtain a restricted Euler-Poisson system (2.19), (2.20), i.e.,

$$\partial_t M + u \cdot \nabla M + M^2 = \frac{k}{n}\rho I_{n\times n},$$
$$\partial_t \rho + u \cdot \nabla \rho + \rho tr M = 0.$$

If we let $\lambda_i(x,t)$ denote the eigenvalues of velocity gradient tensor $\nabla u$, then by the spectral dynamics Lemma 3.1, the eigenvalues and the density $\rho$ are coupled through

$$\partial_t \rho + u \cdot \nabla \rho + \rho \sum_{j=1}^{n} \lambda_j = 0, \tag{7.1}$$

$$\partial_t \lambda_i + u \cdot \nabla \lambda_i + \lambda_i^2 = \frac{k\rho}{n}, \quad i = 1, \cdots, n. \tag{7.2}$$

This is a closed system governing the restricted Euler-Poisson equations, which serves as a simple approximation for the evolution of the full Euler-Poisson system (2.14)-(2.16).

In this section we use the spectral dynamics of the restricted Euler-Poisson equations to show two main points:

1). The global existence of the smooth solutions for a large class of 2D initial configurations– consult Theorem 7.1 below.

2). The finite time blowup of the $n$-dimensional solutions subject to another class of initial data outlined in Theorem 7.4 below. As a consequence of 1) and 2), it follows that the 2D restricted Euler-Poisson equations admit a critical threshold which distinguishes between initial configurations leading to finite time breakdown and global smooth solutions. A detailed study of this 2-dimensional critical threshold phenomena in this context is provided in [29]. This complements the study of critical threshold phenomena for isotropic configurations in the general (global) Euler-Poisson equations presented in [15].

We start with the global regularity of 2-D restricted Euler-Poisson solutions. By well known arguments, the global regularity follows from local existence complemented by a boot-strap argument based on the apriori estimate of $\|\nabla u\|_{L^\infty}$. For the 2D restricted Euler-Poisson model, the velocity gradient tensor $\nabla u$ is completely controlled by its eigenvalues, $\lambda_i, i = 1, 2$, consult [29] for a detailed statement of this argument. With this in mind, it is left to obtain apriori uniform bound of $\lambda_i$'s yielding a sufficient condition for the global existence of smooth solutions for the restricted Euler-Poisson model.

**Theorem 7.1** (Global existence). *The solutions of the 2-D restricted Euler-Poisson equations (7.1),(7.2) remains smooth for all time $t > 0$ if both $\lambda_i(0)$, $i = 1, 2$ are complex, i.e., $Im(\lambda_i(\alpha, 0)) \neq 0$.*

*Proof.* In the 2-D case the density equation (7.1) becomes

$$\rho' + \rho(\lambda_1 + \lambda_2) = 0, \quad ' := \partial_t + u \cdot \nabla_x.$$

¿From (7.2) it follows that the evolution of the divergence $d = \lambda_1 + \lambda_2$, is governed by

$$d' + d^2 - 2\lambda_1\lambda_2 = k\rho,$$

and the evolution of $\Lambda = \lambda_1 \lambda_2$ is given by

$$\Lambda' + d\Lambda = \frac{k}{2}\rho d. \tag{7.3}$$

Introduce the 'indicator' function,

$$\Gamma(t) := \exp(\int_0^t d(x(\alpha, \tau), \tau) d\tau),$$

then the density equation gives

$$\rho(x,t) = \rho_0(\alpha)/\Gamma(\alpha,t), \quad t > 0.$$

Noting that

$$\Gamma' = d\Gamma, \qquad \Gamma'' = (d^2 + d')\Gamma,$$

we then have

(7.4) $$\Gamma'' - 2\Lambda\Gamma = k\rho_0.$$

Substitution of $d = \Gamma'/\Gamma$ and $\rho = \rho_0/\Gamma$ into the $\Lambda$ equation (7.3) it follows that

$$(\Gamma\Lambda)' = \frac{k\rho_0}{2}[ln\Gamma]'.$$

Integration once gives

$$\Gamma\Lambda = \frac{k\rho_0}{2}ln\Gamma + \Lambda_0, \quad \Lambda_0 := (\lambda_1\lambda_2)(\alpha,0),$$

which when inserted into (7.4) yields

$$\Gamma'' = k\rho_0 ln\Gamma + 2\Lambda_0 + k\rho_0.$$

The integral energy becomes

$$[\Gamma']^2 = d_0^2 + 2(2\Lambda_0 + k\rho_0)(\Gamma - 1) + 2k\rho_0 \int_1^\Gamma ln\xi d\xi$$
$$= d_0^2 + 4\Lambda_0(\Gamma - 1) + 2k\rho_0 \Gamma ln\Gamma.$$

Assume that the solution break down at a finite time $t^*$, i.e., $\Gamma(t^*) = 0$, then at this time one has

$$[\Gamma']^2 = d_0^2 - 4\Lambda_0 = (\lambda_{10} - \lambda_{20})^2.$$

Therefore finite time breakdown can not occur if $\lambda_1(\alpha, 0)$ is complex. □

*Remark* 7.2. The above sufficient condition is satisfied, for example, by the initial velocity with large enough vorticity $\omega := u_x - v_y$, associated with the scaled velocity $(u_0(\beta x, y), v_0(x, \beta y))$ with $(u_0(\downarrow, \cdot), v_0(\cdot, \uparrow))$ such that $|\omega_0| \sim \beta^2 >> |d_0|$, implying $Im(\lambda_{10}) \neq 0$.

*Remark* 7.3. What happens with a possible blow-up if both $\lambda_{i0}$ are real? Let $t^*$ be a finite blow-up time. If follows that the blow-up rate is necessarily of the form $(t^* - t)^{-1}$, i.e.,

$$d(t) \sim -\frac{1}{t^* - t} \quad as \quad t \uparrow t^*.$$

This follows from a simple analysis on the following relations

$$0 > \Gamma'(t^*) = \lambda_{10} - \lambda_{20} = \lim_{t \to t^{*-}} d(t)\Gamma(t), \quad \Gamma(t) = \exp(\int_0^t d(\tau)d\tau).$$

The loss of smoothness of the velocity field is closely related with the intricate problem of weak convergence in the absence of the strong convergence. The open question in this context is how the nonlocal term affects the topology of the flow.

To gain further insight on the question of global regularity vs. finite time breakdown, we continue with the $n$-dimensional restricted Euler-Poisson dynamics (7.1)-(7.2). As before, we subtract two consecutive eigenvalue equations in (7.2) to obtain

$$[ln(\lambda_i - \lambda_j)]' + \lambda_i + \lambda_j = 0, \quad \text{for} \quad i \neq j.$$

Summation over $(i, j) \in \mathcal{I}$, with $\mathcal{I}$ defined in (6.8), gives

$$[ln(\amalg_{i(,j) \in \mathcal{I}}(\lambda_i - \lambda_j))]' + N \sum_{k=1}^{n} \lambda_k = 0,$$

Combined with the density equation, $[ln\rho]' + \sum_{k=1}^{n} \lambda_k = 0$, this yields the following global invariants

(7.5) $$\frac{\amalg_{(i,j) \in \mathcal{I}}(\lambda_i - \lambda_j)}{\rho^N} = \text{Const.}$$

The level set of the above invariants is not compact and the finite time singularity can not be ruled out, and in fact, noncompactness implies that certain portion of the phase space must lead to finite time breakdown.

In order to perform a singularity analysis similar to the one provided in §6, we consider a truncated system

$$\rho' = -\rho \sum_{i=1}^{n} \lambda_i,$$

$$\lambda_i' = -\lambda_i^2 \quad i = 1, \cdots, n.$$

Finding its dominant solution of the form

$$(\rho, \lambda_1, \cdots, \lambda_n) \sim (\omega_0, \omega_1, \cdots, \omega_n)\tau^p$$

with $p = (p_0, p_1, \cdots, p_n)$ and $\tau = t^* - t$ leads to

$$-\omega_0 p_0 \tau^{p_0 - 1} = -\omega_0 \tau^{p_0} \sum_{j=1}^{n} \omega_j \tau^{p_j},$$

$$-\omega_i p_i \tau^{p_i - 1} = -\omega_i^2 \tau^{2p_i}.$$

This gives the balance $(\omega, p)$ with

$$\omega = (1, -1, \cdots, -1), \quad p = (-n, -1, \cdots, -1),$$

where $\omega_0 > 0$ is chosen so that it is consistent with the positivity of the density.

Therefore, there exists a general Puiseux-series solution [17] based on the above balance pair, and the blow-up may occur on the orthant $\{+, -, \cdots, -\}$. This, combined with the noncompact integrals derived in (7.5), shows that the solution must exhibit finite-time blow-up in the above orthant.

To summarize, we state the following.

**Theorem 7.4** (Global invariants for $n \geq 2$). *Consider the restricted Euler-Poisson dynamics (7.1)-(7.2) with real initial data $(\rho_0, \lambda_1(0), \cdots, \lambda_n(0))$. Then, there exist $\left[\frac{n}{2}\right]$ global invariants*

$$\text{(7.6)} \qquad \frac{\amalg_{(i,j)\in\mathcal{I}}(\lambda_i - \lambda_j)}{\rho^N} = \text{Const.}, \quad N = \begin{cases} 1, & n \text{ even,} \\ 2, & n \text{ odd.} \end{cases}$$

*The general solution may break down at finite time in the orthant $\{+, -, \cdots, -\}$.*

Two particular cases are worth mentioning. In the 2-D case we have one global invariant $(\lambda_2 - \lambda_1)/\rho = \text{Const.}$, while the global invariant in the 3-D case, corresponding to $\mathcal{I} = \{(1,2), (2,3), (3,1)\}$, is given by

$$\frac{(\lambda_1 - \lambda_2)(\lambda_2 - \lambda_3)(\lambda_3 - \lambda_1)}{\rho^2} = \text{Const.}$$

## 8. Concluding Remarks

This work provides a general framework for several variants of the restricted Euler-dynamics in multi-dimensional case, extending the previous study initiated in [34]. The main tool in this paper is the spectral dynamics analysis. We should point out that this analysis enables us to derive global invariants which are otherwise difficult to detect—one such example was used with the viscous dusty medium model in §5. In particular, we obtain a family of global spectral invariants, interesting for their own sake, for both restricted Euler equations (6.9)-(6.10) and the restricted Euler-Poisson equations (7.6).

Noncompactness of the level set of these global invariants implies the finite time breakdown for a class of initial configurations, for which the local topology of the restricted flow is analyzed. This was demonstrated in Theorem 6.5 in the context of the restricted Euler equations. The finite time breakdown in this restricted model does not necessarily bear on the full, non-restricted Euler equations. On the other extreme we have the possible scenario of a global existence of smooth solutions for restricted models such as restricted Euler-Poisson equations, for which we have the global existence once a critical threshold condition is met. Here, we believe, the global existence does carry over to the question of global existence for the full non-restricted Euler-Poisson equations. In particular, in §7 we have shown the existence of a critical threshold for the 2D restricted model, which in turn suggests the critical threshold phenomena for the full 2D Euler-Poisson equations.

We close this section with following comments.

*Remark* 8.1. Suggesting other nonlocal restricted models. To gain further insight on the fine structure of the flow we propose the following restricted nonlocal models for both Euler-Poisson equations and the incompressible Euler-equations, the analysis of which will appear elsewhere.

- Euler-Poisson equations

    We take the diagonal part of the right side of the $M$-equation in the Euler-Poisson dynamics (2.17), (2.18) to obtain

$$\text{(8.1)} \qquad \partial_t \rho + u \cdot \nabla \rho + \rho \, tr M = 0,$$

$$\text{(8.2)} \qquad \partial_t M + u \cdot \nabla M + M^2 = k(R_i R_j(\rho) \delta_{ij}).$$

- Restricted Euler-dynamics
    A restricted nonlocal Euler dynamics

$$(8.3) \qquad \partial_t M + u \cdot \nabla M + M^2 = (R_i R_j (tr M^2) \delta_{ij}).$$

We note that of course $\partial_t tr M = 0$ and the incompressibility is still invariant.

*Remark* 8.2. Is the spectral dynamics sufficient? We are aware that the spectral dynamics does not tell the whole story for general fluid flows. The following example of a Burgers shear-layer [19, 30] demonstrates this point. Here the simplest solutions of the inviscid Euler equations are the Burgers shear-layer solutions with the velocity field given by

$$u = (h(x_2, t), 0, 0) + (0, -\gamma x_2, \gamma x_3).$$

The velocity gradient tensor,

$$\nabla u = \begin{pmatrix} 0 & h_{x_2} & 0 \\ 0 & -\gamma & 0 \\ 0 & 0 & \gamma \end{pmatrix},$$

has eigenvalues $\{-\gamma, 0, \gamma\}$ which reflect strain effects, but otherwise are independent of the arbitrary shear-layer effect $h(x_2, t)$. Thus, the eigenvalues can not capture the complete behavior of this $h$−dependent flow.

*Remark* 8.3. A main issue in this context is how the restricted Euler-type dynamics related to real flows and at what scale of motion it might apply. In the Navier-Stokes equations, for instance, the nonlocal term should not be ignored at both large and small scales. At large scales the pressure-driven eddy intersections are important and at small scales the velocity gradients are limited by viscous diffusion. We refer the reader to [12] for a detailed discussion on this issue. Another interesting issue left for future research is the recovery of the gradient velocity tensor from the known spectral dynamics.

## 9. Appendix. Trace dynamics for the restricted Euler-equations

This appendix is devoted to an alternative formulation of the spectral dynamics in terms of the traces of $M^k$, $k = 1, \cdots, n$, where $M$ solves the restricted Euler equation

$$(9.1) \qquad \frac{d}{dt} M + M^2 = \frac{1}{n} tr M^2 I_{n \times n}.$$

This is motivated by the trace dynamics originally studied in [34] for $n = 3$.

Here we seek an extension for the general $n$-dimensional setting, which is summarized in the following

**Lemma 9.1.** *Consider the n-dimensional restricted Euler system (9.1) subject to the incompressibility condition $m_1 := tr M = 0$. Then the traces $m_k := tr M^k$ for $k = 2, \cdots, n$ satisfy a closed dynamical system, see (9.2)-(9.4) with (9.6) below, which governs the local topology of the restricted flow.*

*Proof.* Based on the equation (9.1) the transport equations for higher products of $M$ can be written as

$$\frac{d}{dt}M^2 + 2M^3 = \frac{2}{n}MtrM^2,$$

$$\frac{d}{dt}M^3 + 3M^4 = \frac{3}{n}M^2trM^2,$$

$$\cdots$$

$$\frac{d}{dt}M^n + nM^{n+1} = M^{n-1}trM^2.$$

Taking the trace of the above equations and using $m_2 = trM^2$ with $m_1 = 0$ leads to

(9.2) $$\frac{d}{dt}m_2 + 2m_3 = 0,$$

(9.3) $$\frac{d}{dt}m_3 + 3m_4 = \frac{3}{n}m_2^2,$$

$$\cdots$$

(9.4) $$\frac{d}{dt}m_n + nm_{n+1} = m_{n-1}m_2.$$

To close the system, it remains to express $m_{n+1}$ in terms of $(m_1, \cdots, m_n)$. To this end we utilize the Cayley-Hamilton theorem

(9.5) $$M^n + q_1 M^{n-1} + \cdots q_{n-1}M + q_n I = 0,$$

expressed in terms of the characteristic coefficients

$$q_1 = -m_1 = 0, \quad q_2 = -\frac{1}{2}m_2, \quad q_3 = -m_3/3, \quad q_4 = -m_4/4 + m_2^2/8, \quad \cdots.$$

Note that the $q$'s can be expressed in terms of $(m_1, \cdots, m_n)$. (The procedure for computing these coefficients is given at the end of this appendix). Using the Cayley-Hamilton relation (9.5) one may reduce $m_{n+1}$ in (9.4) to lower-order products. In fact, $tr(M \times (9.5))$ gives

(9.6) $$m_{n+1} + q_2 m_{n-1} + \cdots + q_{n-1}m_2 = 0.$$

Substitution into (9.4) yields the closed system we sought for. □

We now turn to consider two examples which demonstrate the above procedure.

**Example 1.** (3-dimensional case $n = 3$, see [34, 5])
In the three dimensional case one has

$$q_1 = 0, \quad q_2 = -\frac{1}{2}m_2, \quad q_3 = det(M) = -\frac{1}{3}m_3,$$

hence

$$M^3 - \frac{1}{2}m_2 M - \frac{1}{3}m_3 = 0.$$

This gives

$$m_4 = \frac{1}{2}m_2^2.$$

Thus a closed system is obtained,

(9.7) $$\frac{d}{dt}m_2 + 2m_3 = 0,$$

(9.8) $$\frac{d}{dt}m_3 + \frac{1}{2}m_2^2 = 0.$$

The invariant of $6m_3^2 = m_2^3 +$Const., could be easily obtained. We consider the phase plane $(m_2, m_3)$, except for the separatrix $6m_3^2 = m_2^3$, all other solutions would not approach the origin and have the finite time breakdown, see Figure 9.1.

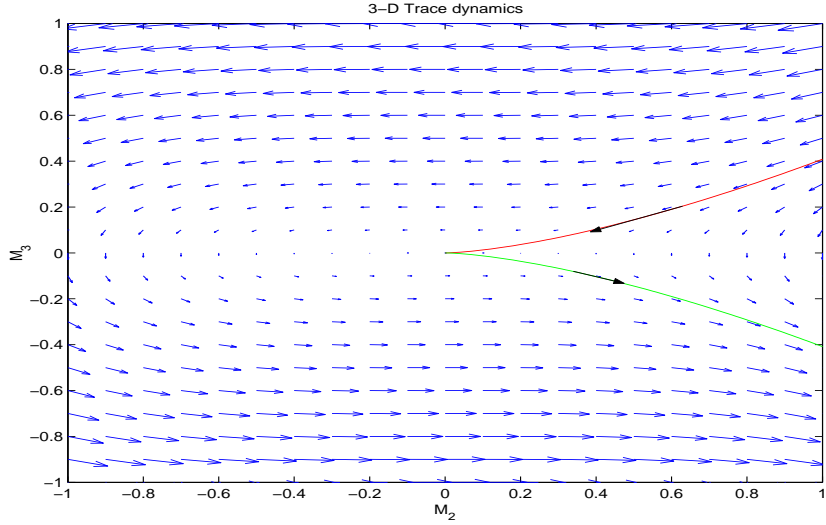

FIGURE 9.1. 3-D Trace-dynamics in Restricted Euler Equations

**Example 2.** (4-dimensional case)

In the four dimensional case one has

$$q_1 = 0, \quad q_2 = -\frac{1}{2}m_2, \quad q_3 = -\frac{1}{3}m_3, \quad q_4 = -\frac{m_4}{4} + \frac{m_2^2}{8}.$$

Hence

$$M^4 - \frac{1}{2}m_2 M^2 - \frac{1}{3}m_3 M - \frac{m_4}{4} + \frac{m_2^2}{8} = 0.$$

Multiplying by $M$ and taking the trace we have

$$m_5 = \frac{1}{2}m_2 m_3 + \frac{1}{3}m_3 m_2 = \frac{5}{6}m_2 m_3.$$

Therefore the resulting closed system becomes

(9.9) $$\frac{d}{dt}m_2 + 2m_3 = 0,$$

(9.10) $$\frac{d}{dt}m_3 + 3m_4 = \frac{3}{4}m_2^2,$$

(9.11) $$\frac{d}{dt}m_4 = -\frac{7}{3}m_3 m_2.$$

This system is still integrable with the following two invariants

$$3m_3^2 = m_2^3 + \frac{3C_1}{4}m_2 + C_2, \quad 12m_4 = 7m_2^2 + C_1,$$

where $C_1, C_2$ are constants integrals of the flow.

*Remark* 9.2. Note that when $C_1 = 0$, the projection of the trajectory on the $m_2 - m_3$ plane has the same topology as that in the 3-D case. See Figure 9.2 for the vector field in $(m_2, m_3, m_4)$ space.

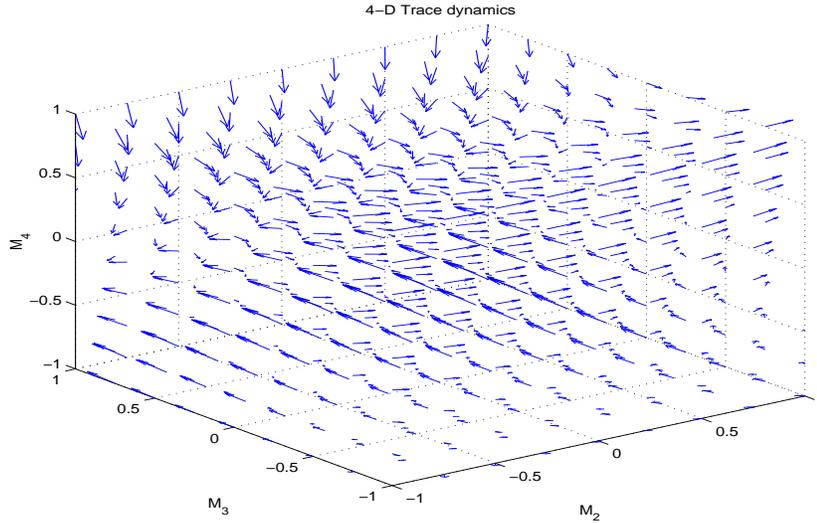

FIGURE 9.2. 4-D Trace-dynamics in Restricted Euler Equations

To gain further insight on the formation of the singularity in this case, we try the dominant solution of the form $\alpha \tau^p$ with $\tau = t^* - t$, for the truncated system

$$\frac{d}{dt}m_2 = -2m_3,$$
$$\frac{d}{dt}m_3 = \frac{3}{4}m_2^2,$$
$$\frac{d}{dt}m_4 = -\frac{7}{3}m_3 m_2.$$

A simple computation gives

$$p = (-2, -3, -4), \quad \alpha = (4, -4, \frac{28}{3}),$$

which shows that the flow may diverge in the orthant $\{+, -, +\}$.

*Remark* 9.3. The above examples demonstrate the difficulty in deriving the global invariants for arbitrary $n > 3$ equations, without the insight provided by the spectral dynamics.

We now conclude this appendix by presenting a procedure of computing the coefficients in the characteristic polynomial for a given matrix.

**Lemma 9.4.** *Let A be a square matrix of order n, its characteristic polynomial reads*

$$det(\lambda I - A) = \sum_{k=0}^{n} q_{n-k}\lambda^k.$$

*Then $q_j = tr(\Lambda^j(A))$, where $\Lambda^j(A)$ is the $j^{th}$ tensor product of A.*

*Proof.* Noting that for $\epsilon = -\lambda^{-1}$

$$det(I + \epsilon A) = \sum_{k=0}^{n} (-1)^k q_k \epsilon^k.$$

On the other hand, if $|\epsilon|$ small, then

$$trlog(I + \epsilon A) = \sum_{k=1}^{\infty} (-1)^{k+1} \frac{tr A^k}{k} \epsilon^k$$

converges. These two relations when combined with the identity

$$det(I + \epsilon A) = \exp\left(trlog(I + \epsilon A)\right)$$

yield

$$\sum_{k=0}^{n} (-1)^k q_k \epsilon^k = \exp\left(\sum_{k=1}^{\infty} (-1)^{k+1} \frac{tr A^k}{k} \epsilon^k\right).$$

Equating the same powers of $\epsilon$ on both sides gives

$$q_0 = 1,$$
$$q_1 = -a_1,$$
$$q_2 = -\frac{a_2}{2} + \frac{a_1^2}{2},$$
$$q_3 = -\frac{a_3}{3} + \frac{a_1 a_2}{2} - \frac{a_1^3}{3!},$$
$$q_4 = -\frac{a_4}{4} + \frac{a_1 a_3}{3} + \frac{a_2^2}{8} - \frac{a_1^2 a_2}{4} + \frac{a_1^4}{4!} \cdots,$$

where $a_k = tr A^k$. This procedure gives the expression of each $q_j$ in terms of $a_k$ for $k = 1, \cdots, n$. □

## ACKNOWLEDGMENTS

Research was supported in part by ONR Grant No. N00014-91-J-1076 (ET) and by NSF grant #DMS01-07917 (ET, HL). We thank Noga Alon for enlightening discussion on counting the different combinatorial arrangements of global invariants in §6.2.

UCLA, Mathematics Department, Los Angeles, CA 90095-1555.
*E-mail address*: `hliu@math.ucla.edu`

UCLA, Mathematics Department, Los Angeles, CA 90095-1555.
*E-mail address*: `tadmor@math.ucla.edu`

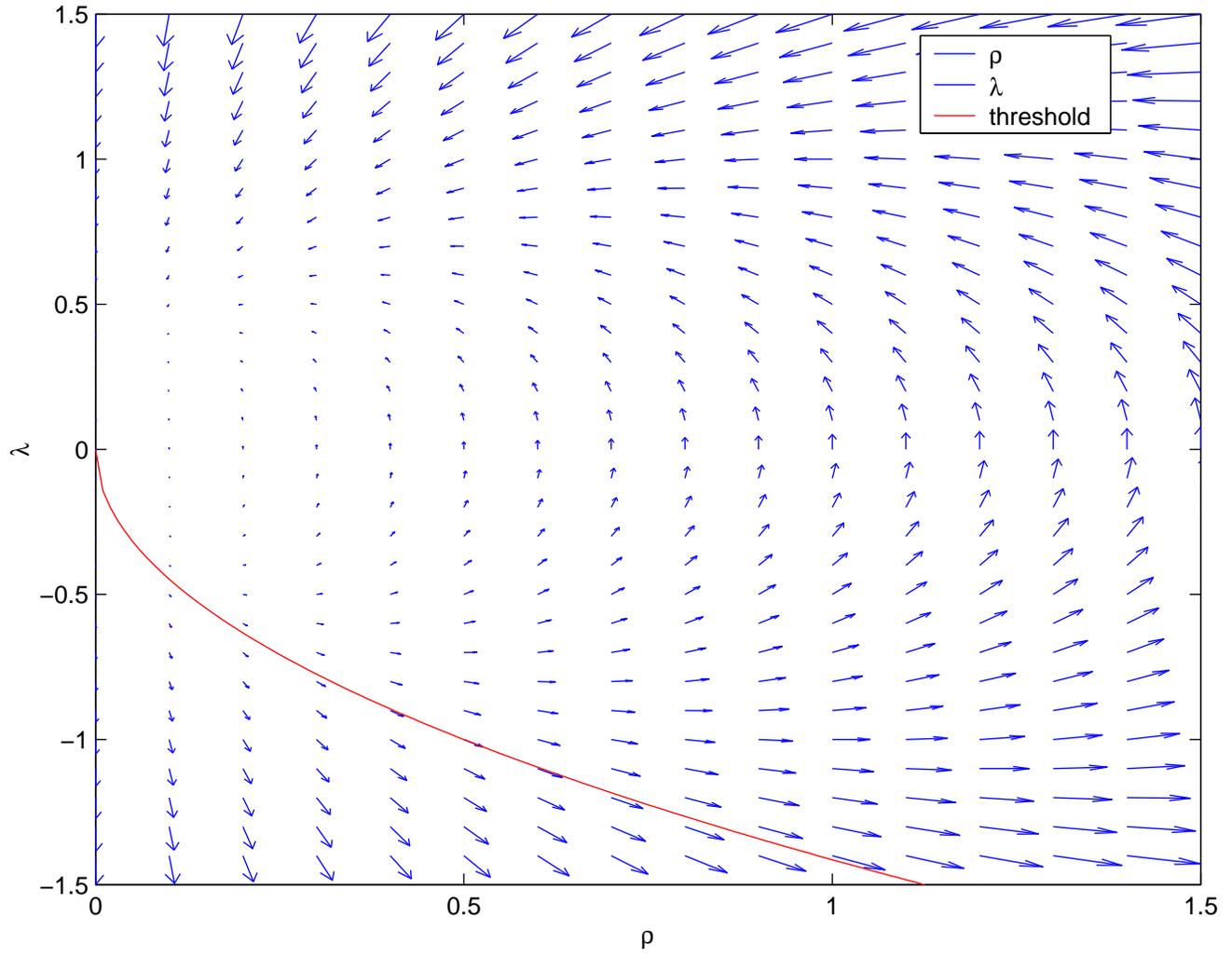